\documentclass[11pt]{article}
\usepackage{enumerate}
\usepackage{amssymb,a4wide,latexsym,makeidx,epsfig,fleqn}
\usepackage{amsthm}
\usepackage{amsmath}
\usepackage{enumerate}
\usepackage{graphicx}
\usepackage{float}
\usepackage{tikz}
\usepackage{hyperref}
\usepackage{blkarray}
\usetikzlibrary{positioning}
\usetikzlibrary{decorations.pathreplacing}
\usepackage{caption}
\allowdisplaybreaks[4]
\newtheorem{theorem}{Theorem}[section]

\newtheorem{lemma}{Lemma}[section]

\usepackage{bm}

\newcommand{\defk}{\operatorname{def}_k}
\newcommand{\odd}{\operatorname{odd}}
\begin{document}
\textwidth 150mm \textheight 225mm
\title{Distance spectral radius conditions for perfect $k$-matching, generalized factor-criticality (bicriticality) and $k$-$d$-criticality of
graphs \thanks{Supported by National Natural Science Foundation of China (No. 12271439).} }
\author{{Kexin Yang$^{a,b}$, Ligong Wang$^{a,b}$\thanks{Corresponding author.
E-mail address: lgwangmath@163.com}, Zhenhao Zhang$^{a,b}$}\\
{\small $^a$School of Mathematics and Statistics, Northwestern
Polytechnical University,}\\ {\small  Xi'an, Shaanxi 710129,
P.R. China}\\
{\small $^b$ Xi'an-Budapest Joint Research Center for Combinatorics, Northwestern
Polytechnical University,}\\
{\small Xi'an, Shaanxi 710129,
P.R. China}\\
{\small E-mail: yangkexi@mail.nwpu.edu.cn, lgwangmath@163.com, zhangzhenhao@mail.nwpu.edu.cn} }
\date{}
\maketitle
\begin{center}
\begin{minipage}{120mm}
\vskip 0.3cm
\begin{center}
{\small {\bf Abstract}}
\end{center}
{\small
	  Let $G$ be a simple connected graph with vertex set $V(G)$ and edge set $E(G)$. A $k$-matching of a graph $G$ is a function
$f:E(G)\rightarrow \{0,1,\ldots, k\}$ satisfying
	  $\sum_{e \in E_G(v)} f(e) \leq k$ for every vertex $v \in V(G)$, where $E_G(v)$ is the set of edges incident with $v$ in $G$. A
$k$-matching of a graph $G$ is perfect if
	  $
	  \sum_{e \in E_G(v) } f(e) = k
	  $
	  for any vertex $v \in V(G)$. The $k$-Berge-Tutte-formula of a graph $G$ is defined as:
	 \[
	 \defk(G) = \max_{S \subseteq V(G)}
	 \begin{cases}
	 	k \cdot i(G - S) - k|S|, & k \text{ is even;} \\[6pt]
	 	\odd(G - S) + k \cdot i(G - S) - k|S|, & k \text{ is odd.}
	 \end{cases}
	 \]
	  A $k$-barrier of the graph $G$ is the subset $S \subseteq V(G)$ that reaches the maximum value in $k$-Berge-Tutte-formula. A connected
graph \( G \) of odd (even) order is a {generalized factor-critical (generalized bicritical) graph about integer \( k \)-matching},
abbreviated as a \( \mathrm{GFC}_k (\mathrm{GBC}_k)\) graph, if $\emptyset$ is a unique $k$-barrier. When $k$ is odd, let \( 1 \leq d \leq k \)
and \( |V(G)| \equiv d \pmod{2} \). If for any \( v \in V(G) \), there exists a \( k \)-matching \( h \) such that
	 $\sum_{e \in E_G(v)} h(e) = k - d$ {and} $\sum_{e \in E_G(u)} h(e) = k$ for any \( u \in V(G) - \{v\} \), then \( G \) is said to be \( k
\)-\( d \)-critical.
	  In this paper, we provide sufficient conditions in terms of  distance spectral radius to ensure that a graph has a perfect $k$-matching
and a graph is \( k \)-\( d \)-critical, $\mathrm{GFC}_k$ or $\mathrm{GBC}_k$, respectively.

}
\end{minipage}
\end{center}

\section{Introduction}
\label{introduction}

All graphs considered in this paper are undirected, connected and simple. Let $G$ be a graph with vertex set $V(G)=\{v_{1}, v_{2}, \ldots,
v_{n}\}$ and edge set $E(G)=\{e_{1},e_{2}, \ldots, e_{m}\}$. We use $n=|V(G)|$ and $m=|E(G)|$ to denote the order and the size of the graph
$G$, respectively. Let $P_n$, $C_n$, $K_n$ and $K_{t, n-t}$ be the path, the cycle, the complete graph and the complete bipartite graph of
order $n$ respectively. Let $G[S]$ be the induced subgraph of $G$ whose vertex set is $S$ and edge set is the edges which have both
endpoints in $S$. For any \( S \subseteq V(G) \), we denote by \( G - S \) the graph obtained from \( G \) by deleting all the vertices in \( S
\) together with all edges incident with vertices in \( S \). A component of a graph $G$ is a maximal connected subgraph of $G$.  A component
is called an odd (even) component if the number of vertices in this component is odd (even). Let \( i(G) \), \( \mathit{odd}(G) \) and \( o(G)
\) denote the number of isolated vertices in \( G \), the number of nontrivial odd components of \( G \) and the number of the odd components
of \( G \), respectively. Then \( o(G) = i(G) + \mathit{odd}(G) \).

We denote by \( G \cup H \) the disjoint union of two graphs \( G \) and \( H \), which is the graph with
$V(G \cup H) = V(G) \cup V(H)$ and $E(G \cup H) = E(G) \cup E(H).$
Denote by \( G \lor H \) the join of two graphs \( G \) and \( H \), which is the graph such that
$V(G \lor H) = V(G) \cup V(H)$ and $E(G \lor H) = E(G) \cup E(H) \cup \{ uv : u \in V(G), v \in V(H) \}.$
The graph \( S_{n,k} \) is obtained from a copy of \( K_k \) by adding \( n - k \) vertices,
each of which has neighborhood \( V(K_k) \), i.e.
$S_{n,k} \cong K_k \lor (n - k)K_1.$

 In 1971, Graham and Pollak {\normalfont \cite{graham and pollack}} firstly studied the distance eigenvalues of graphs. Their research revealed
 that the addressing problem in data communication systems is related to the number of negative distance eigenvalues of trees.
 The distance matrix $D(G)$ of a connected graph $G$ with vertex set $V(G)=\{v_{1}, v_{2}, \ldots, v_{n}\}$ is an $n\times n$  real symmetric matrix
 whose $(i,j)$-entry is $d_G(v_i,v_j)$, where $d_G(v_i,v_j)$ is the distance between the vertices $v_i$ and $v_j$ in the graph $G$. Then we can
 order the eigenvalues of \( D(G) \) as
$
\lambda_1(D(G)) \geq \lambda_2(D(G)) \geq \cdots \geq \lambda_n(D(G))
$. Since the distance matrix $D(G)$ is non-negative and irreducible, by the Perron-Frobenius theorem, \( \lambda_1(D(G)) \) is always positive
(unless \( G \) is trivial) and \( \lambda_1(D(G)) \geq |\lambda_i(D(G))| \) for \( i = 2, 3, \ldots, n \), and we call \( \lambda_1(D(G)) \)
the distance spectral radius of $G$. Furthermore, there exists a unique positive unit eigenvector $\boldsymbol{x}= (x_1,x_2,...,x_n) ^T $
corresponding to \( \lambda_1(D(G)) \), which is called the Perron vector of $D(G)$.

A matching in  a graph \( G \) is a set of pairwise nonadjacent edges of $G$. A perfect matching of $G$ is a matching  which covers every
vertex of \( G \).
In 2005, Brouwer and Haemers {\normalfont \cite{brouwer and haemers}} gave  several spectral sufficient conditions for the existence of a
perfect matching in a graph. After that, the relation between the existence of a perfect matching of a graph and its eigenvalues have received
considerable attention. O {\normalfont \cite{o}} gave a sufficient condition for the existence of perfect matching in terms of the adjacency
spectral radius of a graph $G$. Zhao, Huang and Wang {\normalfont \cite{zhao} generalized this result to the \( A_{\alpha} \)-spectral radius. Zhang and
Lin gave {\normalfont \cite{bi jiao} a distance spectral radius condition to guarantee the existence of a perfect matching in a graph.
		
In 2014, Lu and Wang {\normalfont \cite{lu and wang}} generalized  the concept of perfect matching to perfect $k$-matching.
 A \( k \)-matching of a graph \( G \) is a function \( f: E(G) \to \{0, 1, 2, \ldots, k\} \) satisfying	\(\sum_{e \in E_G(v)} f(e) \leq k\)
 for every vertex \( v \in V(G) \), where \(E_G(v)\) is the set of edges incident with \( v \) in $G$. A \( k \)-matching of \( G \) is perfect
 if \(\sum_{e \in E_G(v)} f(e) = k\) for any vertex \( v \in V(G) \). Clearly, a \( k \)-matching of \( G \) is perfect if and only if 	
 $\sum_{e\in E(G)} f(e) = \frac{k|V(G)|}{2}$.

 When \( k \) is even, Lu and Wang {\normalfont \cite{lu and wang}} proved that a graph has a perfect \( k \)-matching if and only if \( i(G -
 S) \leq |S| \) for all subsets \( S \subseteq V(G) \). A fractional perfect matching of \( G \) is a function \( f : E(G) \to [0, 1] \) such
 that \(\sum_{e \in E_G(v)} f(e) = 1\) for every vetrex \( v \in V(G) \). Scheinerman and Ullman {\normalfont \cite{s and u}} proved that a
 graph \( G \) has a fractional perfect matching if and only if \( i(G - S) \leq |S| \) for all subsets \( S \subseteq V(G) \). Using this
 result, Zhang, Hou and Ren{\normalfont \cite{zhang hou and ren}} gave a distance spectral radius condition to guarantee the existence of a
 fractional perfect matching in a graph. Moreover, Fan, Lin and Lu {\normalfont \cite{fan lin and lu}} combined the adjacency spectral radius
 and minimum degree to provide sufficient conditions for ensuring that a connected graph has a fractional perfect matching. Besides, a graph
 has a perfect \( k \)-matching if and only if it has a fractional perfect matching, where \( k \) is even.

  When \( k \) is odd, Lu and Wang {\normalfont \cite{lu and wang}} also provided a sufficient and necessary condition to ensure the existence
  of perfect \( k \)-matching in a graph, see Lemma \ref{k-matching} in Section 2. Using this result, Zhang and Fan {\normalfont \cite{zhang
  and fan}} provided lower bounds on the number of edges and the adjacency spectral radius for the existence of perfect \( k \)-matching of a
  graph $G$. Niu, Zhang and Wang {\normalfont \cite{niu zhang and wang}}  provided a \( A_{\alpha} \)-spectral radius condition to ensure that a graph
  has a perfect \( k \)-matching, which generalizes the results of Zhang and Fan {\normalfont \cite{zhang and fan}}. Based on the above
  results, we will provide a distance spectral radius condition to ensure the existence of a perfect \( k \)-matching in a graph.

 \noindent\begin{theorem}\label{theorem1}
 	Let $k\geq 1$ be an odd integer. Let $G$ be a connected graph with order $n$ and $n\geq 6$ be an even integer. Then the following
 statements hold.
 	
 	\begin{enumerate}
 		\item[(i)] For $6\le n \leq 8$, if $\lambda_1(D(G)) \le \lambda_1(D(S_{n,\frac{n}{2}-1}))$, then $G$ has a perfect $k$-matching
 unless $G \cong S_{n,\frac{n}{2}-1}$.
 		\item[(ii)] For $n \geq 10$, if $\lambda_1(D(G)) \le \lambda_1(D(G^*))$, then $G$ has a perfect $k$-matching unless $G \cong G^*\cong K_1
 \vee (K_{n-3} \cup 2K_1)$.
 	\end{enumerate}
 \end{theorem}

In order to strengthen the understanding of graph structure with perfect matching, Gallai {\normalfont \cite{gallai}}
and Lov\'{a}sz {\normalfont \cite{lovasz}} introduced the concepts of factor-critical and bicritical graphs, respectively. A graph \( G \) is called to be {factor-critical} if \( G - v \) has a perfect matching for every \( v \in V(G) \). A graph \( G \) is {bicritical} if \( G - x -
y \) has a perfect matching for any two distinct vertices \( x \) and \( y \) in \( G \). Favaron {\normalfont \cite{favaron and s}} extended
them to \( k \)-factor-critical in 1986. A graph \( G \) is {\( k \)-factor-critical} if \( G - S \) has a perfect matching for every subset \(
S \subseteq V(G) \) with \( |S| = k \). The study of the relationship between the  $k$-factor-critical and the spectral radius has received
widespread attention, as we can see \cite{fan and lin, niu zhang and wang,zhang and fan,zeng}.
In 2026, Chang, Liu and Qiu \cite{chang liu and qiu} generalized the definitions and results about factor-critical graphs and bicritical graphs
to the corresponding graphs about integer $k$-matching respectively. They also gave the definition of the $k$-$d$-critical graph, where $k\geq
3$ is odd.  Before introducing the relevant definitions, we firstly present the \( k \)-Berge--Tutte formula. Liu and Liu \cite{liu1,liu2} gave
the \( k \)-Berge--Tutte formula of a graph \( G \) as follows.

\[
\defk(G) = \max_{S \subseteq V(G)}
\begin{cases}
	k \cdot i(G - S) - k|S|, & k \text{ is even;} \\[6pt]
	\odd(G - S) + k \cdot i(G - S) - k|S|, & k \text{ is odd.}
\end{cases}
\]

The vertex subset \( S \subseteq V(G) \) that reaches the maximum value in the \( k \)-Berge--Tutte formula is called a $k$-barrier of a graph
$G$.  A connected graph \( G \) of odd order is a {generalized factor-critical graph about integer \( k \)-matching}, abbreviated as a \(
\mathrm{GFC}_k \) graph, if \( G \) has no non-empty \( k \)-barrier. A connected  graph \( G \) of even order is a {generalized bicritical
graph about integer \( k \)-matching}, abbreviated as a \( \mathrm{GBC}_k \) graph, if \( G \) has no non-empty \( k \)-barrier.

When \( k\geq 3 \) is  an odd integer, Liu, Su and Xiong \cite{liu2} showed the following lemma.

\begin{lemma}[{\cite{liu2}}]\label{prop:1.1}
	Let \( G \) be a connected odd graph with at least three vertices. Then \( G \) is \( \mathrm{GFC}_k \) if and only if for any \( v \in
V(G) \), there exists an integer \( k \)-matching \( h \) such that
	\[\sum_{e \in E_G(v)} h(e) = k - 1 \quad \text{and} \quad \sum_{e \in E_G(u)} h(e) = k \quad \text{for any other vertex } u.\]
\end{lemma}

By Lemma~\ref{prop:1.1}, Chang, Liu and Qiu \cite{chang liu and qiu} extended the above results and defined a \( k \)-\( d \)-critical graph.

For an odd integer \( k \geq 3 \), let \( 1 \leq d \leq k \) and \( |V(G)| \equiv d \pmod{2} \). If for any \( v \in V(G) \), there exists a \(
k \)-matching \( h \) such that
\[
\sum_{e \in E_G(v)} h(e) = k - d \quad \text{and} \quad \sum_{e \in E_G(u)} h(e) = k
\]
for any \( u \in V(G) - \{v\} \), then \( G \) is said to be \( k \)-\( d \)-critical.

Chang, Liu and Qiu \cite{chang liu and qiu} provided the necessary and sufficient conditions for a graph to be $\mathrm{GFC}_k$,
$\mathrm{GBC}_k$ or \( k \)-\( d \)-critical, respectively, see Lemma \ref{generalized k-factor} and Lemma \ref{kd} in Section 2. Zhang and
Wang \cite{zhang wang} established tight sufficient conditions in terms of size and adjacency spectral radius for a graph \( G \) to be \( k
\)-\( d \)-critical, $\mathrm{GFC}_k$ or $\mathrm{GBC}_k$, respectively.
 Motivated from the results, we will give sufficient conditions in terms of the distance spectral
radius to ensure that a graph is \( k \)-\( d \)-critical, $\mathrm{GFC}_k$ or $\mathrm{GBC}_k$, respectively. However, if \( d = k \), the
complete graph \( K_n \) is not \( k \)-\( d \)-critical and the complete graph attains the minimum distance spectral radius among all graphs with $n$ vertices.
Then we consider \( 1 \leq d < k \) as follows.

Firstly, when $k\geq 3$ is odd, we present the distance spectral radius condition for a graph to be \( k \)-\( d \)-critical.

\noindent\begin{theorem}\label{t5}
	Let $k\geq 3$ be an odd integer. Let $G$ be a connected graph with order $n\geq 3$. For \( 1 \leq d < k \) and \( n \equiv d \pmod{2} \),
then the following statements hold.
	\begin{enumerate}
		\item[(i)] If $ \lambda_1(D(G))\le \lambda_1(D(K_1 \vee (K_{n-2} \cup K_1)))$, then \( G \) is \( k \)-\( d \)-critical for odd
$n\geq 3$ and even $n\geq 10$ unless $G \cong K_1 \vee (K_{n-2} \cup  K_1)$.
		\item[(ii)] If $\lambda_1(D(G))\le\lambda_1(D( S_{n,\frac{n}{2}})))$, then \( G \) is \( k \)-\( d \)-critical for even $4 \le n\le
8$ unless $G \cong S_{n,\frac{n}{2}}$.
	\end{enumerate}
\end{theorem}

When \( k\geq 3 \) is odd, Chang, Liu and Qiu \cite{chang liu and qiu} also proved that for any \( k \)-\( d \)-critical graph \( G \), when \(
d \) is odd, \( G \) is \( \mathrm{GFC}_k \); when \( d \) is even, \( G \) is \( \mathrm{GBC}_k \). Therefore, let $d=1$ and $d=2$ in Theorem
\ref{t5} to obtain the corresponding results for a graph to be \( \mathrm{GFC}_k \) or \( \mathrm{GBC}_k \), see Theorems \ref{theorem3} and
\ref{theorem4}.

\noindent\begin{theorem}\label{theorem3}
	Let both $k\geq 3$ and $n\geq 3$ be odd integers. Let $G$ be a connected graph with order $n$. If $ \lambda_1(D(G))\le \lambda_1(D(K_1 \vee
(K_{n-2} \cup K_1)))$, then \( G \) is $ \mathrm{GFC}_k$ unless $G \cong K_1 \vee (K_{n-2} \cup  K_1)$.
	
\end{theorem}

 \noindent\begin{theorem}\label{theorem4}
 	Let $k\geq 3$ be an odd integer. Let $G$ be a connected graph with order $n$ and $n\geq 4$ be an even integer. Then the following
 statements hold.
 	\begin{enumerate}
 	\item[(i)] For $n\geq 10$, if $ \lambda_1(D(G))\le \lambda_1(D(K_1 \vee (K_{n-2} \cup K_1)))$, then \( G \) is  $\mathrm{GBC}_k$ unless
 $G \cong K_1 \vee (K_{n-2} \cup  K_1)$.
 	\item[(ii)] For $ 4 \le n\le 8$, if $\lambda_1(D(G))\le\lambda_1(D( S_{n,\frac{n}{2}})))$, then \( G \) is $\mathrm{GBC}_k$ unless $G
 \cong S_{n,\frac{n}{2}}$.
 	\end{enumerate}
 \end{theorem}
Secondly, when  \( k\geq 2 \) is even, we present the distance spectral radius condition for a graph to be  $\mathrm{GFC}_k$ or
$\mathrm{GBC}_k$.

\noindent\begin{theorem}\label{theorem2}
	Let $k\geq 2$ be an even integer. Let $G$ be a connected graph with order $n\geq 3$. Then the following statements hold.
	\begin{enumerate}
		\item[(i)] If $ \lambda_1(D(G))\le \lambda_1(D(K_1 \vee (K_{n-2} \cup K_1)))$, then \( G \) is $ \mathrm{GFC}_k$ for odd $n\geq 3$
and $G$ is $\mathrm{GBC}_k$ for even $n\geq 10$ unless $G \cong K_1 \vee (K_{n-2} \cup  K_1)$.
		\item[(ii)] If $\lambda_1(D(G))\le\lambda_1(D( S_{n,\frac{n}{2}})))$, then \( G \) is \( \mathrm{GBC}_k \) for even $4 \le n\le 8$
unless $G \cong S_{n,\frac{n}{2}}$.
	\end{enumerate}
\end{theorem}

 The rest of this paper is organized as follows. In Section 2, we present some important preliminary results that will be used later. In
 Section 3,
 we prove Theorem \ref{theorem1}. In Section 4, we prove Theorem \ref{t5}. then let $d=1$ and $d=2$, we obtain the corresponding
 results (Theorem \ref{theorem3} and \ref{theorem4}) for a graph to be \( \mathrm{GFC}_k \) and \( \mathrm{GBC}_k \), respectively. In Section 5,
 we prove Theorem \ref{theorem2}.

\section{Preliminaries}
\label{sec:Preliminaries}

In this section, we will present some important results that will be used in our subsequent arguments.

Let \( M \) be a real \( n \times n \) matrix, and let \( X = \{1, 2, \ldots, n\} \). Given a partition \(\Pi : X = X_1 \cup X_2 \cup \cdots
\cup X_k \), the matrix \( M \) can be correspondingly partitioned as

\[
M =
\begin{pmatrix}
	M_{1,1} & M_{1,2} & \cdots & M_{1,k} \\
	M_{2,1} & M_{2,2} & \cdots & M_{2,k} \\
	\vdots & \vdots & \ddots & \vdots \\
	M_{k,1} & M_{k,2} & \cdots & M_{k,k}
\end{pmatrix}
\]

The quotient matrix of \( M \) with respect to \(\Pi\) is defined as the \( k \times k \) matrix \( B_{\Pi} = (b_{i,j})_{k \times k} \) where
\( b_{i,j} \) is the average value of all row sums of \( M_{i,j} \). If each block \( M_{i,j} \) of \( M \) has constant row sum \( b_{i,j} \),
then we call the partition \(\Pi\) an equitable partition, and \( B_{\Pi} \) the equitable quotient matrix of $M$.

\noindent\begin{lemma}\label{quotient}{\normalfont(\cite{quotient})}   Let \( M \) be a real symmetric matrix, and let \( \lambda_1(M) \) be
the largest eigenvalue of \( M \). If \( B_{\Pi} \) is an equitable quotient matrix of \( M \), then the eigenvalues of \( B_{\Pi} \) are also
eigenvalues of \( M \). Furthermore, if \( M \) is nonnegative and irreducible, then \( \lambda_1(M) = \lambda_1(B_{\Pi}) \).
\end{lemma}

We present a fundamental result to compare the distance spectral radii of a graph and its spanning subgraph,
which is a corollary of the Perron-Frobenius theorem.
\noindent\begin{lemma}\label{bian duo d da}{\normalfont(\cite{bian duo d da})}
	Let $e$ be an edge of a graph $G$ such that $G - e$ is connected. Then
	\[
	\lambda_1(D(G)) < \lambda_1(D(G - e)).
	\]
\end{lemma}

The following lemma gives a sufficient and necessary condition to ensure that a graph has a perfect $k$-matching.
\noindent\begin{lemma}\label{k-matching}{\normalfont(\cite{lu and wang})}
Let $k > 1$ be an odd integer. A graph $G$ contains a perfect $k$-matching if and only if
\[
\operatorname{odd}(G - S) + k \cdot i(G - S) \le k\cdot|S|
\]
for all subsets $S \subseteq V(G)$.
\end{lemma}

Chang, Liu and Qiu {\normalfont \cite{chang liu and qiu}} obtained the following sufficient and necessary conditions for a graph to be a
generalized factor-critical graph, a generalized bicritical graph and a \( k \)-\( d \)-critical graph, respectively.

\noindent\begin{lemma}\label{generalized k-factor}{(\cite{chang liu and qiu})}
	Let \( G \) be a connected graph with order $n\geq 3$ and \( k \geq 2 \). Then
	
	(i) When \( k \) is even, \( G \) is \( \mathrm{GFC}_k \, (\mathrm{GBC}_k) \) if and only if $n$ is odd (even) and \( i(G - S) \leq |S| - 1
\) for any \( \emptyset \neq S \subset V(G) \).
	
	(ii) When \( k \) is odd, \( G \) is \( \mathrm{GFC}_k \) if and only if $n$ is odd and \( \mathrm{odd}(G - S) + k \cdot i(G - S) \leq k|S|
- 1 \) for any \( \emptyset \neq S \subset V(G) \).
	
	(iii) When \( k \) is odd, \( G \) is \( \mathrm{GBC}_k \) if and only if $n$ is even and \( \mathrm{odd}(G - S) + k \cdot i(G - S) \leq
k|S| - 2 \) for any \( \emptyset \neq S \subset V(G) \).
	
\end{lemma}

\noindent\begin{lemma}\label{kd}{(\cite{chang liu and qiu})}
	Let \( G \) be a graph of order \( n \geq 3 \) and \( k \geq 3 \) be odd. For \( 1 \leq d \leq k \) and \( n \equiv d \pmod{2} \), then \(
G \) is \( k \)-\( d \)-critical if and only if
	\[
	\operatorname{odd}(G - S) + k \cdot i(G - S) \leq k|S| - d \quad \text{for any } \emptyset \neq S \subseteq V(G).
	\]
	
\end{lemma}

Let $W(G) = \sum_{i<j} d_{G}(v_i,v_j)$ be the Wiener index of a connected graph $G$ of order $n$ with vertex set $V(G)=\{v_{1}, v_{2}, \ldots,
v_{n}\}$. Note that $\lambda_1(D(G)) = \max_{x \in \mathbb{R}^n} \frac{x^T D(G) x}{x^T x}$.

Then we have
$$
\lambda_1(D(G)) = \max_{x \in \mathbb{R}^n} \frac{x^T D(G) x}{x^T x} \geq \frac{\bm{1}^T D(G) \bm{1}}{\bm{1}^T \bm{1}} \geq \frac{2W(G)}{n},
$$
where $\bm{1} = (1, 1, \dots, 1)^T$.

\noindent\begin{lemma}\label{bi jiao}(\cite{bi jiao}, pp. 317--319) Let \( p \geq 2 \) and \( n_i \geq 1 \) for \( i = 1, \ldots, p \). If
\(\sum_{i=1}^p n_i = n - s \) where \( s \geq 1 \), then
	\[\lambda_1 (D(K_s \vee (K_{n-s-p+1} \cup (p-1)K_1)))\le
	\lambda_1 (D(K_s \vee (K_{n_1} \cup K_{n_2} \cup \cdots \cup K_{n_p})))
	\]
	with equality if and only if \( n_i = 1 \) for \( i = 2, \ldots, p \).
\end{lemma}

\noindent\begin{lemma}\label{perfect}(\cite{bi jiao}, pp. 319--320) Let \( s \geq 1 \) be an integer and $G$ be a connected graph of even order
$n$ and $n\geq 2s+2$. Then the following statements hold.
	\begin{enumerate}
		\item[(i)] If $n \geq 2s+4$, then $\lambda_1(D(K_1 \vee (K_{n-3} \cup 2K_1))) \le \lambda_1(D( K_s \vee \left( K_{n-2s-1} \cup (s +
1)K_1 \right)))$ with equality if and only if $s=1$.
		\item[(ii)] If $n=2s+2$ and $4\le n\le 8$, then $\lambda_1(D( S_{n,\frac{n}{2}-1} )))\le \lambda_1(D(K_1 \vee (K_{n-3} \cup 2K_1)))$
with equality if and only if $n=4$.
		\item[(iii)] If $n=2s+2$ and $n\geq 10$, then $\lambda_1(D(K_1 \vee (K_{n-3} \cup 2K_1)))< \lambda_1(D( S_{n,\frac{n}{2}-1})))$.
		
	\end{enumerate}
\end{lemma}

To prove our main theorems, we first prove the following lemma.

\noindent\begin{lemma}\label{bi jiao 2} Let \( s \geq 1 \) be an integer and $G$ be a connected graph of order $n$ and $n\geq 2s$. Then the
following statements hold.
	\begin{enumerate}
		\item[(i)] If $n \geq 2s+2$, then $\lambda_1(D(K_1 \vee (K_{n-2} \cup K_1))) \le \lambda_1(D( K_s \vee \left( K_{n-2s} \cup s K_1
\right)))$ with equality if and only if $s=1$.
		\item[(ii)] If $n=2s+1$, then $\lambda_1(D(K_1 \vee (K_{n-2} \cup K_1))) \le \lambda_1(D( S_{n,\frac{n-1}{2}})))$  with equality if
and only if $n=3$.
		\item[(iii)] If $n=2s$ and $n\geq 10$, then $\lambda_1(D(K_1 \vee (K_{n-2} \cup K_1))) <\lambda_1(D( S_{n,\frac{n}{2}})))$.
		\item[(iv)] If $n=2s$ and $2\le n\le 8$, then $\lambda_1(D( S_{n,\frac{n}{2}})))\le \lambda_1(D(K_1 \vee (K_{n-2} \cup K_1))) $ with
equality if and only if $n=2$.
	\end{enumerate}
\end{lemma}
\begin{Tproof}\textbf{.}
	
	Let \( G^s= K_s \lor (K_{n-2s} \cup sK_1) \) and $G^*= K_1 \lor (K_{n-2} \cup K_1)$.
	We will  compare the distance spectral radii of $G^s$ and $G^*$ as follows.
	
\textbf{Claim 1.} If $n \geq 2s+2$, then $\lambda_1(D(G^*)) \le \lambda_1(D(G^s))$ with equality if and only if $s=1$.

 If \(s = 1\), then \({G^*} \cong G^s\). Now we suppose \(s\geq 2\), so that $n\geq 2s \geq 4$. Next we write the quotient matrix of the
 partition \(\{V(K_{n-2s}), V(K_{s}), V(sK_1)\}\) of \(G^s\) is
	
	\[
	\begin{pmatrix}
		n-2s-1 & s & 2s \\
		n-2s & s-1 & s \\
		2(n-2s) & s & 2(s-1)
	\end{pmatrix}.
	\]
	
	The characteristic polynomial of the above matrix is
	\begin{align*}
		f(x) &=x^3- (n+s-4)x^2+ (5-n (2 s+3)+s (5 s-1))+n ((s-2) s-2)x+s^2 (5-2 s)+2.
	\end{align*}

	We easily obtain that the largest root of \(f(x) = 0\) is \(\lambda_1(D({G^s}))\). Let $s=1$, since \(\lambda_1(D(G^*)) = \theta(n)\)
(simply \(\theta\)) is the largest root of the equation $q(x)=x^3+(3-n)x^2+(9-5 n)x-3 n+5$, we have
	\begin{align*}
		h(\theta) &=f(\theta)-q(\theta)  \\
		&= (s-1) \left(-\theta^2-(2 n-5 s-4)\theta+n s-n-2 s^2+3 s+3\right).
	\end{align*}
	
	Moreover, \(\theta = \lambda_1(D(G^*)) \geq \frac{2W(G^*)}{n} = \frac{n^2+n-4}{n} \geq n \) and $s \geq 2$. Next we only need to prove
\(h_1(\theta) = -\theta^2-(2 n-5 s-4)\theta+n s-n-2 s^2+3 s+3<0 \) when \(\theta \geq n \), then \(h(\theta) < h_1(\theta) < 0\). Then we can
obtain \(\lambda_1(D({G^s})) \geq \lambda_1(D(G^*))\). Note that
	$$\frac{-2n+5s+4}{2}=-n+\frac{5s}{2}+2<n.$$
	
So when \(\theta \geq n \), $h_1(\theta)$ is monotonically decreasing as $\theta$ increases and $h_1(\theta)\le h_1(n)$. Next we only need to
prove $h_1(n)<0$.
	
Define $g(n)=h_1(n)=-3 n^2+ (6 s+3)n+(3-2 s)s+3$. Then
		$$\frac{6s+3}{6}=s+\frac{1}{2}<2s+2 .$$
So when $n\geq 2s+2$, $g(n)$ is monotonically decreasing and $h(\theta)<h_1(\theta)\le g(2s+2)=-(2 s+3)s-3<0$. So Claim 1 holds.

\textbf{Claim 2.}  If $n = 2s+1$, then $\lambda_1(D(G^*)) \le \lambda_1(D(S_{n, \frac{n-1}{2}}))$ with equality if and only if $n=3$.
	
Let \( n = 2s + 1 \), then \( g(n) = -\frac{n^2}{2}+\frac{5 n}{2}+1 < 0 \) when \( n \geq 6 \). If \( n = 2s + 1 \), we have \( G^s\cong S_{n,
\frac{n-1}{2}} \). Then we write the quotient matrix of the partition $\{V(K_{\frac{n-1}{2}}), V((\frac{n+1}{2})K_1)\}$
of \( S_{n, \frac{n-1}{2}} \) is
\[\begin{pmatrix}
	\frac{n-3}{2} & \frac{n+1}{2} \\
	\frac{n-1}{2}  & n-1
\end{pmatrix}.\]

By a simple calculation,

\[\lambda_1(D(S_{n, \frac{n-1}{2}})) = \frac{1}{4} \left(\sqrt{5 n^2+2 n-3}+3 n-5\right).\]	

 Besides, we know that \( \lambda_1(D(G^*)) = \theta(n) \) is the largest root of the equation $q(x)=x^3+(3-n)x^2+(9-5 n)x-3 n+5$. If \( s = 1
 \) and \( n = 3 \), then \( G^* \cong S_{n, \frac{n-1}{2}} \). If \( s = 2 \) and \( n=5\), then \( q(x) = x^3-2 x^2-16 x-10 \) and \(
 \lambda_1(D(S_{5,2})) = \frac{2 \sqrt{33}+10}{4} \). It is obvious that \( q (\frac{2 \sqrt{33}+10}{4})  = \frac{ \sqrt{33}-3}{2} > 0 \) and
 \( q'\left(\frac{2 \sqrt{33}+10}{4}\right) = \frac{11 \sqrt{33}+35}{2} > 0 \), so \( \lambda_1(D(S_{5,2})) > \theta(5) \). So Claim 2 holds.
	
\textbf{Claim 3.} If \( n = 2s \) and \( n \geq 10 \), then $\lambda_1(D(G^*)) < \lambda_1(D(S_{n, \frac{n}{2}}))$.
	
Let \( n = 2s \), then \( g(n) = -\frac{n^2}{2}+\frac{9 n}{2}+3 < 0 \) when \( n \geq 10 \).  So Claim 3 holds.

\textbf{Claim 4.} If \( n = 2s \) and \(2 \le n \le 8 \), then $\lambda_1(D(G^*)) \geq \lambda_1(D(S_{n, \frac{n}{2}}))$ with equality if and
only if $n=2$.

If \( s = 1 \) and \( n = 2 \), then \( G^* \cong S_{n, \frac{n}{2}} \). If \( n = 2s \), we have \( G^s \cong S_{n, \frac{n}{2}} \). Then we
write the quotient matrix of the partition $\{V(K_{\frac{n}{2}}), V((\frac{n}{2})K_1)\}$
of \( S_{n, \frac{n}{2}} \) is
\[\begin{pmatrix}
	\frac{n}{2}-1 & \frac{n}{2} \\
	\frac{n}{2}  & n-2
\end{pmatrix}.\]

By a simple calculation,

\[\lambda_1(D(S_{n, \frac{n}{2}})) = \frac{1}{4} \left(\sqrt{5 n^2-4 n+4}+3 n-6\right).\]	 Besides, we know that \( \lambda_1(D(G^*)) =
\theta(n) \) is the largest root of the equation $q(x)=x^3+(3-n)x^2+(9-5 n)x-3 n+5$.
If \( s = 2 \) and \( n=4\), then \( q(x) =x^3-x^2-11 x-7 \) and \( \lambda_1(D(S_{2,2})) = \frac{2 \sqrt{17}+6}{4} \). It is obvious that \( q
(\frac{2 \sqrt{17}+6}{4})  = \frac{-3(\sqrt{17}+5)}{2} < 0 \), so \( \lambda_1(D(S_{2,2})) < \theta(4) \). If \( s = 3 \) and \( n=6\), then \(
q(x) =x^3-3 x^2-21 x-13 \) and \( \lambda_1(D(S_{3,3})) =  \sqrt{10}+3 \). It is obvious that \( q (\sqrt{10}+3)  = -2 \left(\sqrt{10}+8\right)
< 0 \), so \( \lambda_1(D(S_{3,3})) < \theta(4) \).
	If \( s = 4 \) and \( n=8\), then \( q(x) =x^3-5 x^2-31 x-19 \) and \( \lambda_1(D(S_{4,4})) = \frac{2 \sqrt{73}+18}{4} \). It is obvious
that \( q (\frac{2 \sqrt{73}+18}{4})  = \frac{3\left(\sqrt{73}-9\right)}{2} < 0 \), so \( \lambda_1(D(S_{4,4})) < \theta(8) \). So Claim 4
holds.  \hfill $\square$
	\end{Tproof}

\section{Proof of Theorem \ref{theorem1}}

In this section, we give the proof of Theorem \ref{theorem1}, which provides a sufficient condition for a graph to have a perfect $k$-matching
in terms of the distance spectral radius of a graph.

\begin{Tproof}\textbf{.}
	By way of contradiction assume that a graph \( G \) of even order $n$ contains no perfect \( k \)-matching with the minimum distance
spectral radius. By Lemma \ref{k-matching}, there exists a subset \( S \subseteq V(G) \) such that
	\[
	\operatorname{odd}(G - S) + k \cdot i(G - S) \geq k \cdot |S|+1.
	\]
	
	 Let \( p = \operatorname{odd}(G - S) \), \( i = i(G - S) \) and \( s = |S| \), then \( p + ki \geq ks+1 \). Note that \( G \) is
connected, which implies that $S$ is not empty and $s\geq 1$. Since \( n \) is even and \( k \) is odd, then \( p + ki \) and \( ks \) have the
same parity. Therefore,
	\[
	p + ki \geq ks + 2.
	\]
	
By Lemma \ref{bian duo d da}, we claim that the induced graph \( G[S] \) is complete, each component of \( G - S \) is either a complete graph
or an isolated vertex, and each vertex of \( S \) is adjacent to each vertex of \( G - S \). Next, we consider the following two cases.

\noindent\textbf{Case 1.} $	\operatorname{odd}(G - S)=0$.

In this case, \( p = \operatorname{odd}(G - S) = 0 \), then  we obtain \( ki \geq ks + 2 \) in the graph $G$. Thus \( i \geq s + 1 \). Since \(
n \) is even, then \( i \) and \( s \) have the same parity. We have \( i \geq s + 2 \) in $G$. Next we claim that \( G - S \) has no even
component. By way of contradiction assume that there is an even component in \( G - S \). Let \( G' \) be a graph obtained from \( G \) by
adding an edge between the even component and a trivial component of \( G - S \). In this case, we have
	\[
	\operatorname{odd}(G' - S) + k \cdot i(G' - S)=1+ k \cdot(i(G - S)-1) \geq ks + 2.
	\]

By Lemma \ref{k-matching}, \( G' \) has no perfect \( k \)-matching. By Lemma \ref{bian duo d da}, we have \( \lambda_{1}(D(G')) <
\lambda_{1}(D(G)) \), which contradicts the choice of \( G \). Then \( G \cong G_1 \cong K_s \vee (n - s)K_1 \).

\noindent\textbf{Case 2.} $\operatorname{odd}(G - S)\neq 0$.

In this case, \( p = \operatorname{odd}(G - S) \geq 1 \). Next we claim that \( G - S \) has no even component. By way of contradiction assume
that there is an even component in \( G - S \). Let \( G'' \) be a graph obtained from \( G \) by adding an edge between the even
component and a non-trivial odd component of \( G - S \).  then we obtain
\[
\operatorname{odd}(G'' - S) + k \cdot i(G'' - S)=\operatorname{odd}(G - S)+ k \cdot i(G - S) \geq ks + 2.
\]

By Lemma \ref{k-matching}, \( G'' \) has no perfect \( k \)-matching. By Lemma \ref{bian duo d da}, we have \( \lambda_{1}(D(G'')) <
\lambda_{1}(D(G)) \), a contradiction.

Then we obtain that each component of \( G - S \) is odd. It is clear that $$G =G_2= K_s \vee \left( K_{n_1} \cup K_{n_2} \cup \cdots \cup
K_{n_p} \cup iK_1 \right),$$
where $n_1 \geq n_2 \geq \cdots \geq n_p\geq 1$, $i = i(G - S)$ and \( s = |S| \).

Let \( G_3 = K_s \vee \left( K_{n-s-t} \cup tK_1 \right) \), where \( t = i + p - 1 \). Note that \( p + ki \geq ks + 2 \). We obtain
	\[
\operatorname{odd}(G_3 - S) + k \cdot i(G_3 - S)=1 + kt = 1 + k(i + p - 1) = p + ki + (k - 1)(p - 1) \geq ks + 2.
	\]
	
By Lemma \ref{k-matching}, \( G_3 \) contains no perfect \( k \)-matching. By Lemma \ref{bi jiao}, \( \lambda_{1}(D(G_3)) \leq
\lambda_{1}(D(G_2)) \), then $G=G_3$.

Let \( G^s = K_s \vee \left( K_{n-2s-1} \cup (s + 1)K_1 \right) \).  Since \( 1 + kt \geq ks + 2 \) in the graph $G_3$, we have \( t \geq s + 1
\); since \( k(n - s) \geq ks + 2 \) in the graph $G_1$ and \( s \) and \( n - s \) have the same parity, we have \( n - s \geq s + 2 \). So \(
G_3 \) and $G_1$ are both spanning subgraphs of the graph \( G^s \). By Lemma \ref{bian duo d da}, then \( \lambda_{1}(D(G^s)) \leq
\lambda_{1}(D(G_1)) \) and \( \lambda_{1}(D(G^s)) \leq \lambda_{1}(D(G_3)) \). We obtain $G\cong G^s = K_s \vee \left( K_{n-2s-1} \cup (s +
1)K_1 \right) $.

Next we claim that \( n \geq 2s + 2 \). If \( i \geq s + 1 \), then \( n \geq s + i \geq 2s + 1 \). Since \( n \) is even, we have \( n \geq 2s
+ 2 \). If \( i \leq s \), then
\[
n \geq s + p + i \geq s + k(s - i) + 2 + i = 2s + 2 + (k - 1)(s - i) \geq 2s + 2.
\] 	

Let \( G^* = K_1 \vee (K_{n-3} \cup 2K_1) \).  Note that \( G^* = G^s\) if \( s = 1 \). In the sequel, we compare the size of \(
\lambda_{1}(D(G^*)) \) and \( \lambda_{1}(D(G^s)) \), where \( s \geq 2 \).

 If $n=2s+2$, then $G^s\cong  S_{n,\frac{n}{2}-1}$. Let $S = V(K_{\frac{n}{2}-1})$. Then
	\[
	\operatorname{odd}(S_{n,\frac{n}{2}-1} - S) + k \cdot i(S_{n,\frac{n}{2}-1} - S) = k(\frac{n}{2}+1) \geq k(\frac{n}{2}-1) + 2.
	\]
	
 By Lemma \ref{k-matching} \( S_{n,\frac{n}{2}-1} \) contains no perfect \( k \)-matching. Therefore, by Lemma \ref{perfect}, then Theorem
 \ref{theorem1} holds. 	\hfill $\square$

\end{Tproof}

\section{Proof of Theorem \ref{t5}}
In this section, we give the proof of Theorem \ref{t5}, which provides a sufficient condition in terms of the distance spectral radius for a
graph of order $n\geq 3$ to be \( k \)-\( d \)-critical, where $k\geq 3$ is an odd integer, $1 \leq d < k$ and $n \equiv d \pmod{2}$.
\begin{Tproof}\textbf{.}
	
	By way of contradiction assume that a graph \( G \) of order $n\geq 3$ is not $k$-$d$-critical with the minimum distance spectral radius.
By Lemma \ref{kd}, there exists a subset \( \emptyset \neq S \subset V(G) \) such that
	\[
	\operatorname{odd}(G - S) + k \cdot i(G - S) \geq k |S|-d+1.
	\]
	
	Let \( p = \operatorname{odd}(G - S) \), \( i = i(G - S) \) and \( s = |S| \), then \( p + ki \geq ks-d+1 \). Note that $S \neq \emptyset$,
then $s\geq 1$. According to the assumption, we have $k\geq 3$ is an odd integer, $1 \leq d < k$ and $n \equiv d \pmod{2}$.
	
	By Lemma \ref{bian duo d da}, we claim that the induced graph \( G[S] \) is complete, each component of \( G - S \) is either a complete
graph or an isolated vertex, and each vertex of \( S \) is adjacent to each vertex of \( G - S \). Next, we consider the following two cases.
	
	\noindent\textbf{Case 1.} $\operatorname{odd}(G - S)=0$.
	
	In this case, since $p=\operatorname{odd}(G - S)=0$ and $1 \leq d < k$, then $ki \geq ks-d+1$ and \( i \geq s \geq1 \).
	
	If $i=s$, then  we firstly claim that \( G - S \) has at most one even component. Suppose that \( G-S \) has two or more even components.
Let \( G' \) be a graph obtained from $G$ by adding edges which connect all even components. Then we have
	\[
	\operatorname{odd}(G' - S) + k \cdot i(G' - S)=k \cdot i(G - S) \geq ks-d+1.
	\]
	
	By Lemma \ref{kd}, \( G' \) is not $k$-$d$-critical. By Lemma \ref{bian duo d da}, we have \( \lambda_{1}(D(G')) < \lambda_{1}(D(G)) \),
which contradicts the choice of \( G \). Then $ G \cong G_1= K_s \vee ( K_{n-2s} \cup sK_1)$.

	If $i\geq s+1$, we firstly claim that \( G - S \) has no even component. By way of contradiction assume that there is an even component in
\( G - S \). Let \( G'' \) be a graph obtained from \( G \) by adding an edge between the even component and an isolated vertex  of \( G
- S \). Then we have
	\[
	\operatorname{odd}(G'' - S) + k \cdot i(G'' - S)=1+ k \cdot(i(G - S)-1) \geq 1+ks\geq ks-d+1.
	\]
	
	By Lemma \ref{kd}, \( G'' \) is not $k$-$d$-critical. By Lemma \ref{bian duo d da}, we have \( \lambda_{1}(D(G'')) < \lambda_{1}(D(G)) \),
which contradicts the choice of \( G \). Then \( G \cong G_2 \cong K_s \vee (n - s)K_1 \).
	
	Since $G_2$ is a spanning subgraph of the graph $G_1$, by Lemma \ref{bian duo d da}, we have \( \lambda_{1}(D(G_1)) < \lambda_{1}(D(G_2))
\). Then we consider $G \cong G_1= K_s \vee ( K_{n-2s} \cup sK_1)$ in this case.

	\noindent\textbf{Case 2.} $	\operatorname{odd}(G - S)\neq 0$.
	
	We firstly claim that \( G - S \) has no even component. By way of contradiction assume that there is an even component in \( G - S \). Let
\( G''' \) be a graph obtained from \( G \) by adding an edge between the even component and a non-trivial odd component of \( G - S \).
In this case, \( p = \operatorname{odd}(G - S) \geq 1 \), then we obtain
	\[
	\operatorname{odd}(G''' - S) + k \cdot i(G''' - S)=\operatorname{odd}(G - S)+ k \cdot i(G - S) \geq ks-d+1.
	\]
	
	By Lemma \ref{kd}, \( G''' \) is not $k$-$d$-critical. By Lemma \ref{bian duo d da}, we have \( \lambda_{1}(D(G''')) < \lambda_{1}(D(G))
\), a contradiction. Then we obtain that each component of \( G - S \) is odd. It is clear that $$G =G_3= K_s \vee \left( K_{n_1} \cup K_{n_2}
\cup \cdots \cup K_{n_p} \cup iK_1 \right),$$
	where $n_1 \geq n_2 \geq \cdots \geq n_p\geq 1$, $i = i(G - S)$ and \( s = |S| \).
	
	Let \( G_4 = K_s \vee \left( K_{n-s-t} \cup tK_1 \right) \), where \( t = i + p - 1 \). Note that \( p + ki \geq ks-d+1 \). We obtain
	\[
	\operatorname{odd}(G_4 - S) + k \cdot i(G_4 - S)=1 + kt = 1 + k(i + p - 1) = p + ki + (k - 1)(p - 1) \geq ks-d+1.
	\]
	
	By Lemma \ref{kd}, \( G_4 \) is not $k$-$d$-critical. By Lemma \ref{bi jiao}, \( \lambda_{1}(D(G_4)) \leq \lambda_{1}(D(G_3)) \), then
$G=G_4$.
	
	Let \( G^s =G_1= K_s \vee \left( K_{n-2s} \cup sK_1 \right) \).  Since \( 1 + kt \geq ks-d+1 \) and $1 \le d < k$ in the graph $G_4$, we
have \( t \geq s  \). Then \( G_4\) is a spanning subgraph of the graph \( G^s \). By Lemma \ref{bian duo d da}, Then \( \lambda_{1}(D(G^s))
\leq \lambda_{1}(D(G_4)) \). We obtain $G\cong G^s = K_s \vee \left( K_{n-2s} \cup sK_1 \right) $.
	
	Let \( G^* = K_1 \vee (K_{n-2} \cup K_1) \).  Note that \( G^* = G^s\) if \( s = 1 \). In the sequel, we compare the size of \(
\lambda_{1}(D(G^*)) \) and \( \lambda_{1}(D(G_s)) \), where \( s \geq 2 \).
	Therefore, by Lemma \ref{bi jiao 2}, then Theorem \ref{t5} holds.
	\hfill $\square$
\end{Tproof}

\section{Proof of Theorem \ref{theorem2}}

In this section, we give the proof of Theorem \ref{theorem2}, which provides a sufficient condition in terms of the distance spectral radius
for a graph with odd (even) order $n\geq 3$ to be \( \mathrm{GFC}_k \, (\mathrm{GBC}_k) \), where $k\geq 2$ is an even integer.

\begin{Tproof}\textbf{.}
	By way of contradiction assume that a graph \( G \) of odd (even) order $n$ is \( \mathrm{GFC}_k \, (\mathrm{GBC}_k) \) with the minimum
distance spectral radius. By Lemma \ref{generalized k-factor}, there exists a subset \( \emptyset \neq S \subset V(G) \) such that
	\[
	 i(G - S) \geq  |S|.
	\]
	
	Let \( i = i(G - S) \) and \( s = |S| \). Then \( i \geq s \) and $s\geq 1$. Hence $n\geq i+s \geq 2s$.

\noindent\textbf{Claim 1.} $G-S$ contains at most one non-trivial connected component.
	
	Suppose that \( G-S \) has two or more non-trivial connected components. Let \( G' \) be a graph obtained from $G$ by adding edges which
connect all non-trivial connected components. Then we have $i(G' - S)=i(G - S) \geq |S|$. By Lemma \ref{generalized k-factor}, \( G' \) is not
\( \mathrm{GFC}_k \, (\mathrm{GBC}_k) \).  By Lemma \ref{bian duo d da}, we have \( \lambda_{1}(D(G')) < \lambda_{1}(D(G)) \), which
contradicts the choice of \( G \). So Claim 1 holds.

	By Lemma \ref{bian duo d da}, we claim that the induced graph \( G[S] \) is complete, each component of \( G - S \) is either a complete
graph or an isolated vertex, and each vertex of \( S \) is adjacent to each vertex of \( G - S \). Next, we consider two cases to prove this
theorem.
	
	\noindent\textbf{Case 1.} 	$G-S$ has exactly one non-trivial component.
	
	In this case, we have \( G \cong K_s \vee (K_{n_1} \cup iK_1) \), where \( n_1 = n-s-i \geq 2 \). If \( i \geq s+1 \), we can obtain a new
graph \( G'' \) by adding edges between the non-trivial component and an isolated vertex. Then we obtain the graph \( G'' \) satisfies
\( i(G''-S) = i(G-S)-1 = i-1 \geq s \). By Lemma \ref{generalized k-factor}, \( G'' \) is not \( \mathrm{GFC}_k \, (\mathrm{GBC}_k) \).  By
Lemma \ref{bian duo d da}, we have \( \lambda_{1}(D(G'')) < \lambda_{1}(D(G)) \), which contradicts the choice of \( G \). Hence, we have  \( i
= s \) and so \( G \cong G^s =K_s \vee (K_{n-2s} \cup sK_1) \).

\noindent\textbf{Case 2.} $G-S$ has no non-trivial component.
		
	In this case, we have \( G \cong G_1 = K_s \vee (n - s)K_1 \). Since $s\geq 1$, then  $G_1$ is a spanning subgraph of the  graph $G^s$. By
Lemma \ref{bian duo d da}, then \( \lambda_{1}(D(G^s)) < \lambda_{1}(D(G_1)) \). Therefore, we only consider \( G \cong G^s=K_s \vee (K_{n-2s}
\cup sK_1) \).

Let \( G^* = K_1 \vee (K_{n-2} \cup K_1) \).  Note that \( G^* = G^s \) if \( s = 1 \). In the sequel, we compare the size of \(
\lambda_{1}(D(G^*)) \) and \( \lambda_{1}(D(G^s)) \), where \( s \geq 2 \).

If $n=2s+1$, then $G^s\cong  S_{n,\frac{n-1}{2}}$.	
Let $S = V(K_{\frac{n-1}{2}})$. Then
	\[
 i(S_{n,\frac{n-1}{2}} - S) = \frac{n+1}{2} \geq \frac{n-1}{2}.
	\]
	
	By Lemma \ref{generalized k-factor}, the graph \( S_{n,\frac{n-1}{2}} \) is not $\mathrm{GFC}_k $.
	
	If $n=2s$, then $G^s\cong  S_{n,\frac{n}{2}}$.	
	Let $S = V(K_{\frac{n}{2}})$. Then
	\[
	i(S_{n,\frac{n}{2}} - S) = \frac{n}{2} \geq \frac{n}{2}.
	\]
	
	By Lemma \ref{generalized k-factor}, the graph \( S_{n,\frac{n}{2}} \) is not \(\mathrm{GBC}_k\).
	
	Therefore, by Lemma \ref{bi jiao 2}, then Theorem \ref{theorem2} holds.
	
		\hfill $\square$

\end{Tproof}

\noindent\textbf{Declaration of competing interest}

 There is no competing interest.

\noindent\textbf{Data availability}

No data was used for the research described in the article.

\end{document}